\documentclass[a4paper, 12pt]{article}

\usepackage{graphicx}
\usepackage{amssymb}

\def\blankbox{{ ~\hfill$\rlap{$\sqcap$}\sqcup$}}

\begin{document}


\bigskip \bigskip

\centerline {\bf The Second Subconstituent of some Strongly Regular Graphs}

\bigskip \bigskip

\centerline{Norman Biggs}

\bigskip \bigskip

\centerline{Department of Mathematics}

\centerline{London School of Economics}

\centerline{Houghton Street}

\centerline{London WC2A 2AE}

\centerline{U.K.}

\centerline{n.l.biggs@lse.ac.uk}

\bigskip

\centerline{February 2010}

\bigskip \bigskip

\centerline{\bf Abstract}
\bigskip

This is a report on a failed attempt to construct new graphs $X^h$ that are strongly regular with
parameters $((h^4 + 3h^2 +4)/2, \, h^2+1,\,  0, 2)$.
The approach is based on the assumption that the second subconstituent of $X^h$ has an equitable
partition with four parts.
For infinitely many odd prime power values of $h$ we construct a graph $G^h$ that is a plausible candidate
for the second subconstituent. Unfortunately
 we also show that the corresponding $X^h$ is strongly regular
only when  $h=3$, in which case the graph is already known.

\vfill \eject

 {\bf 1. Introduction}
\medskip

We should like to be able to construct graphs $X$ that have the following properties:
\smallskip

{\leftskip 30pt

$\bullet$  $X$ is regular with degree $k$;

$\bullet$  $X$ is triangle-free;

$\bullet$  any two non-adjacent vertices have just two common neighbours.

\par}
\smallskip

Standard calculations with eigenvalues and multiplicities show that  $k-1$ must be a square $h^2$, with $h$
not congruent to $0$ modulo $4$. In the standard terminology [5,11], $X$ is a strongly regular
graph with parameters
$$((h^4 + 3h^2 +4)/2, \, h^2+1,\,  0, 2).$$
 Although there are infinitely many possibilities, only a few graphs are known, even when the number of common
  neighbours is allowed to be an arbitrary constant $c \neq 2$ [10].
The topic is particularly interesting because the known graphs are associated with remarkable groups.
\smallskip

For each vertex $v$ of $X$ we denote by  $X_1(v)$, $X_2(v)$ the sets of vertices at distance $1,2$ respectively
from $v$.  We call the graph induced by $X_2(v)$ the {\it second subconstituent} of $X$.  We shall usually
write it as $X_2$, although there is no reason why it should be independent of $v$.
The results in [2,11] establish that
 $X_2$ is a connected graph of degree $k-2$ with diameter 2 or 3.
Furthermore, the only numbers that can be eigenvalues of $X_2$ are: $k-2$, $-2$ and the
eigenvalues $\lambda_1$, $\lambda_2$ of $X$.
\smallskip

There are three known examples:  $h=1,2,3$ corresponding to $k=2,5,10$.
When $k=2$ we have the $4$-cycle.  When $k=5$ we have a graph with 16 vertices known as the
Clebsch graph, a name suggested by Coxeter [6] because the graph represents a geometrical
configuration discussed
by Clebsch.  For this graph $X_2$ is the Petersen graph.
When $k=10$ we have the Gewirtz graph with $56$ vertices.  It can be represented by
taking the vertices to be a  set of $56$ ovals in PG$(2,4)$, and making two vertices adjacent when
the corresponding ovals are disjoint. (An historical note about this graph is appended to this
paper.)  The algebraic properties of the Gewirtz graph have been studied in
detail by Brouwer and Haemers [4], and a list of the $56$ ovals may be found at [12].
Using this list, it can be verified that the
45 vertices of $X_2$ are partitioned with respect to
their distance from a given vertex $w \in X_2$ as  $\{w\} \cup P \cup S$, where  $|P| = 8$ and
$|S| = 36$.  The 36 vertices are of two types.
One set $Q$ of 16 vertices has the property that each is adjacent to 1 vertex in $P$,
while the complementary set $R$ of 20 vertices is such that each is adjacent to 2 vertices in $P$.
Further analysis shows that the partition  $\{w\}\cup P \cup Q \cup R$ is equitable [11, p.195],
with the intersection numbers given by the matrix

$$\pmatrix{ 0 &1 &0 &0 \cr
             8 &0 &1 &2 \cr
             0 &2 &2 &4 \cr
             0 &5 &5 &2 \cr }.$$

The eigenvalues of this intersection matrix are $8, 2, -2, -4$. The numbers $2$ and $-4$ are the eigenvalues
of the Gewirtz graph, while $8$ $(=k-2)$ and $-2$ are the only other eigenvalues permitted by the general
 theorem mentioned above.
On this basis, it seems  worthwhile to investigate possible generalizations.
\bigskip

{\bf 2.  Properties of the second subconstituent}
\medskip

We begin by constructing a suitable intersection matrix for the second subconstituent,
for a general value of $k$.

Let $G = (V,E)$ be a graph with vertex-set $V = K^{(2)}$, the set of unordered pairs of
elements of a set $K$,
where $|K| = k$, and suppose the edge-set $E$ is defined so that the following conditions hold.
\smallskip

{\leftskip 20pt

{\bf C1}  {\it For each vertex $ab \in V$ there is a partition of $V$ with four parts,}
$$ V = \{ab\} \cup P_{ab} \cup Q_{ab} \cup R_{ab}$$
{\it such that $P_{ab}  = \{cd \mid \{ab,cd\} \in E\}$ and $Q_{ab} = \{cd \mid |ab \cap cd| = 1\}$.}
\smallskip

{\bf C2}  {\it This partition is equitable with intersection matrix}
$$M = \pmatrix { 0    &1    &0     &0   \cr
             k-2  &0    &1     &2   \cr
             0    &2    &2     &4   \cr
             0    &k-5  &k-5   &k-8 \cr }.$$
\par}

\medskip

The fact that $P_{ab}$ and $Q_{ab}$ are disjoint implies that if $\{ab,cd\}$ is an edge, then
$a,b,c,d$ are distinct. It follows from the definition of $Q_{ab}$ that $|Q_{ab}| = 2(k-2)$.
The other parameters then imply that
$$|P_{ab}| = k-2, \quad |R_{ab}| = \frac{1}{2}(k-2)(k-5).$$
It is easy to see that a graph $G$ with the given properties would be triangle-free,  regular with degree $k-2$,
and have diameter 2.
\smallskip

The conditions are clearly meaningful only when $k \ge 8$. Since we know that the corresponding
strongly regular graphs can exist only when $k=h^2 + 1$, and that they do exist when
$k=2$ and $k=5$,  we shall assume that $k\ge 10$
 in what follows.
\medskip

{\bf Theorem 1} \quad Let $G$ be a graph satisfying conditions {\bf C1} and {\bf C2}, and
let $C$ be the bipartite graph (claw) with vertex-set $\{*\} \cup K$. Denote by $G \oplus C$  the graph formed
from the disjoint union of $G$ and $C$ by adding edges joining each vertex $ab$ in $G$ to the
vertices $a$ and $b$ in $C$.  Then $G \oplus C$ is a strongly regular graph with degree $k$, it is triangle-free,
and each pair of non-adjacent vertices has just 2 common neighbours.
\smallskip

{\it Proof} \quad Since $G$ is regular with degree
$k-2$,  $G \oplus C$ is regular with degree $k$. Since $C$ is a claw, there are no triangles containing
the vertex $*$.
A triangle containing the vertex $a \in K$ would have to contain vertices $ab$ and $ab'$ in $G$,
but these vertices are not adjacent. Finally, a triangle containing
the vertex $ab$ would lie wholly in $G$, but $G$ is triangle-free.
\smallskip

It remains to check that any two non-adjacent vertices in $G \oplus C$ have exactly two common neighbours.
If the two vertices are $*$ and $ab$, the neighbours are $a$ and $b$, and if the two vertices are $a$,
$b$, the common neighbours are $*$ and $ab$. If the two vertices are of the form $ab$ and $ac$, then
$ac \in Q_{ab}$, and the neighbours are $a$ and the unique vertex in $P_{ab}$ that is adjacent to $ac$.
If the two vertices are of the form $ab$ and $cd$, where $cd \in R_{ab}$, then
the neighbours are the two vertices in $P_{ab}$ that are adjacent to $cd$.
\blankbox

\bigskip

{\bf 3. Construction of triangle-free graphs $G^q$}
\medskip

We attempt to construct graphs satisfying conditions {\bf C1} and {\bf C2}. We know that
$k-1$ must be a square, say $h^2$. Let
$h= q$, where $q$ is a prime power, so that  $k= q^2 + 1$.  Take $K$ to be the set of points
on the projective line $PG(1, q^2)$, that is

$$K= {\mathbb F}_{q^2} \cup \{\infty\} = \langle t \rangle \cup \{0, \infty \}.$$

Here ${\mathbb F}_{q^2}$ is the finite field of order $q^2$,  $\infty$ is the
conventional `point at infinity', and  $t$ is a primitive element of the field, so that
$\langle t \rangle$ is a cyclic group of order $q^2 -1 = k-2$.
The group  PGL$(2, q^2)$ of projective linear transformations acts 3-transitively on $K$,
and hence transitively on the
unordered pairs $ab$ in $K^{(2)}$, and the stabilizer  of the pair $0\infty$ is generated by
$x \mapsto tx$ and
$x \mapsto x^{-1}$.  When $q$ is an odd prime power (so that $k$ is even) its orbits on $K^{(2)}$ are as follows:

$$\{0 \infty\}, \quad O_0 = \{0x \mid x \in \langle t \rangle\} \cup
\{ \infty x \mid x \in \langle t \rangle\},$$
and $\frac{1}{2} (k-2)$ orbits of the form

$$O_v = \{ vx\; x \mid x \in \langle t \rangle \} \quad (v = t, t^2, \ldots, t^{(k-2)/2}).$$

The orbit $O_0$ has size $2(k-2)$, the orbits $O_v$ ($v \neq -1$) have size $k-2$, and the orbit
$O_{-1}$ has size $(k-2)/2$. Note that when $q$ is a power of $2$ the orbit-partition takes a slightly
different form: this is consistent with the fact that no construction can work when
$k-1= h^2$ with $h \equiv 0$ (mod $4$), by the feasibility conditions. (The exceptional case $q=2$
has already been covered.)
\smallskip

The idea of the following  construction is to define a graph with vertex-set $V = K^{(2)}$ such
that for a suitable value of $u$, the partition
postulated in condition {\bf C1} (taking the vertex $ab$ to be $0\infty$) is given by

$$P= O_u, \qquad Q = O_0, \qquad  R = \bigcup_{v \neq u, 0} O_v.$$

The construction depends on the {\it cross-ratio}, which is defined for any points $a,b,c,d \in K$ by the rule

$$ (ab|cd)   = \frac{(a-c)(b-d)}{(a-d)(b-c)},$$

with the usual conventions about $\infty$. The cross-ratio is $1$ if and only if
$a=b$ or $c=d$ or both, and so this value does not occur when $ab$ and $cd$ are in
$V = K^{(2)}$.  Given the unordered pair of unordered pairs $ab$ and $cd$,
the cross-ratio $(ab|cd)$ takes only two values $\rho$ and $\rho^{-1}$, which it is convenient to write
in the form $(ab|cd) = \rho^{\pm}$.
\smallskip

Let $V = K^{(2)}$, and given $u \in \langle t \rangle, u \neq \pm 1$ define $E_u$ to be the set of pairs
$\{ab, cd\}$ such that $(ab|cd) = u^{\pm}$.  Since $(0 \infty \mid ux\, x) = u$, it follows that in the
graph $G_u = (V, E_u)$ the set of vertices adjacent to $0\infty$ is the orbit $O_u$, as defined above.
\smallskip

We consider the possibility that, for a suitable value of $u$,  $G_u$ is a graph in which the partition
given above is equitable, with the intersection matrix $M$ as in condition {\bf C2}.
The first step is to ensure that $M_{PP} = 0$, which means that $G_u$ is triangle-free.
\medskip

{\bf Lemma} \quad  Let

$$ \Omega \; = \; \{v \in K \mid v= (x+x^{-1} -1)^{\pm} \; {\rm for\; some} \;
x \in \langle t \rangle, x \neq 1\}. $$

Then the graph $G_u$ is triangle-free if and only if $u$ is not in $\Omega$.
\smallskip

{\it Proof} \quad The group PGL$(2, q^2)$ acts as  a group of automorphisms of $G_u$ since it preserves
cross-ratios, and so $G_u$ is vertex-transitive.  Hence we need only consider the possibility of triangles
containing a given vertex, say $0\infty$.
The stabilizer of $0 \infty$ contains
$x \mapsto tx$ and $x \mapsto 1/x$, and so we can assume that two edges of the triangle are
 $\{ 0 \infty,  u \;1\}$ and $\{ 0 \infty,  ux \;x \}$, with $x \neq 1$.
Now, the vertices $u\;1$ and $ux\;x$ are adjacent if and only if

$$\frac{ (ux-u)(u-1)}{ (ux-1)(x-u) }  = u^{\pm}.$$

After some rearrangement, this reduces to

$$ u = (x + x^{-1} -1)^{\pm}.$$

In other words, there is a triangle if and only if $u$ is in $\Omega$.

\blankbox
\bigskip

Since $x + x^{-1} -1$ is symmetrical with respect to inverting $x$, the set $\Omega$ can be
found by calculating at the values of $x + x^{-1} -1$
for $x= t^j$, $j = 1, 2, 3, \ldots, (k-2)/2$.  For example, when the field is ${\mathbb F}_9$ with
the primitive element $t$ satisfying $t^2 + t + 2 =0$,  we have the table

$$\matrix{ j\quad       &x        &x+ x^{-1} - 1  &(x+x^{-1}-1)^{-1}  \cr
                        &         &           &                 \cr
           1\quad            &t        &2t         &2+2t             \cr
           2\quad            &1+2t     &2          &2                \cr
           3\quad            &2+2t     &1+t        &t                \cr
           4\quad            &2        &0          &\infty           \cr }.$$

Since $1+2t$ and its inverse $2+t$ are not in $\Omega$, we conclude that the graph $G_{1+2t}$ is triangle-free.
\smallskip

Generally, the special orbit $O_u$ could be any one of the $O_v$ except $O_0$ and $O_{-1}$, thus $\{u, u^{-1}\}$
could be any one of the $(k-4)/2$ pairs $\{t^j, t^{-j} \}$, $ j= 1, 2, \ldots,(k-4)/2$.
As a working definition let us say that $u$ is {\it admissible} if (1) $u = \neq 0, \infty, -1, 1$, and
(2) $\{u, u^{-1}\} \notin \Omega$.
At first sight it appears that as many as $(k-2)/2$ pairs are not admissible, because they are in $\Omega$,
but fortunately things are not so bad.
\medskip

{\bf Theorem 2} \quad Suppose that $q$ is an odd prime power and there is an element
$\zeta$ in ${\mathbb F}_{q^2}$ such that $\zeta^2 = 3$.
Then there is at least one $u$ of the form $t^j$ with $j \in \{1,2, \ldots, (k-4)/2\}$ such that
$u$ is not in $\Omega$, and hence the graph $G_u$ is triangle-free.
\smallskip

{\it Proof} \quad When $q$ is odd  $q^2$ is congruent to $1$ mod $4$, and
so there is an element $\iota$ such that $\iota^2 = -1$.  Then

$$\iota + \iota^{-1} - 1 = (\iota + \iota^{-1} -1)^{-1} = -1,$$
which means that the `pair' $\{-1, -1\}$ occurs in $\Omega$. But
 $u = -1$ is not admissible anyway, and so
 the number of non-admissible pairs in $\Omega$ is effectively reduced to at most $(k-4)/2$.
\smallskip

Similarly, if we can find an $x$ such that $x + x^{-1} -1 = 0$  then the pair $\{0, \infty\}$ will
occur in $\Omega$, and since this pair is also not admissible,  the number
of non-admissible pairs in $\Omega$ will be reduced to at most $(k-6)/2$.
It is easy to see that this happens if
there is an element $\zeta \in {\mathbb F}_{q^2}$ such that $\zeta^2 = 3$.  In that case, let

$$\theta = 2^{-1}(1+\iota \zeta), \quad {\rm so\;that} \quad
\theta^{-1} = 2^{-1}(1- \iota\zeta) \quad {\rm and} \quad \theta + \theta^{-1} -1 = 0.$$

Hence an admissible $u$ must exist.
\blankbox
\bigskip

For example, in the field ${\mathbb F}_{25}$ with
the primitive element $t$ satisfying $t^2 + t + 2 =0$, we have $\iota =2$,
$\zeta= 4+3t$, $\theta = 2+3t$.   Hence the pair $\{2+3t, 4+2t\}$ is admissible.
A complete check shows that there are two other
admissible pairs  $\{1+2t, 2+4t\}$, and $\{3+t, 4+3t\}$.
\bigskip

It is easy to see that there are infinitely many fields ${\mathbb F}_{q^2}$ which contain
an element with $\zeta^2 = 3$, for example by applying the law of quadratic reciprocity.
We do not pursue this matter, since the final step is to rule out all fields except ${\mathbb F}_9$, where
explicit calculation shows that the construction works.

\bigskip

{\bf 4.  Failure of the construction in general}
\medskip

We now know that in many cases a graph $G_u$ can be constructed satisfying the condition $M_{PP} = 0$. But
it remains to check that the other entries of $M$ are correct.  In fact, several of them
can be verified, but it turns out that
the condition $M_{PR} =2$ cannot be satisfied in general.
\medskip

{\bf Theorem 3} \quad Let $G_u$ be defined for an odd prime power $q$ as in Section 3, and let the
partition $\{0 \infty\}\cup P \cup Q \cup R$ of the
vertices of $G_u$ be as stated there. Then this partition is equitable with an intersection
matrix of the form required by condition {\bf C2} only when $q =3$.
\smallskip

{\it Proof} \quad We shall show that the condition $M_{PR}=2$ cannot hold, except when $q=3$.
\smallskip

A typical vertex in $R$ is $vw\, w$ where $v \neq u, 0, \infty, 1$ and $w \in \langle t \rangle$ .
This vertex is adjacent to the vertices $ux\, x$ in $P$ for which
$(ux\, x \mid vw\, w) = u$ or $u^{-1}$. These two equations can be written as quadratics in $x$:

$$ux^2 - (1+u) vwx + vw^2 = 0, \qquad ux^2 -(1+u)wx + vw^2 =0,$$

and their discriminants are

$$\Delta^+ = w^2((1+u)^2v^2 - 4uv), \qquad \Delta^- = w^2((1+u)^2 - 4uv).$$

If there are just two solutions for $x$ then either
(1) exactly one of $\Delta^+$, $\Delta^-$ is a square in ${\mathbb F}_{q^2}$, or (2) $\Delta^+$ and
$\Delta^-$ are both zero.
\smallskip

Consider first the case $v = -1$.  Here $\Delta^+ = \Delta^-$ and so their common value must be zero.
That is,
$$(1 + u^2) +4u =0.$$

Then for all $v \neq -1$
$$\Delta^+ = w^2 (-4uv^2 -4uv) = v w^2 (-4u - 4uv) = v \Delta^-.$$

Hence in order that exactly one of $\Delta^+$, $\Delta^-$ is a square, $v$ must be a non-square, and this must
hold for all the orbits $O_v \subseteq R$ except $O_{-1}$. In the case $q=3$ we chose $u = 1+t = t^2$, so
$R = O_t \cup O_{t^3} \cup O_{-1}$, and the condition is satisfied.  But for $q \ge 5$ there must be
at least one square among the relevant values of $v$ and the condition cannot be satisfied.
\blankbox

\bigskip

{\bf Historical note on the Gewirtz graph}
\medskip

Gewirtz discussed his graph in two papers published in 1969 [8,9].
Brouwer [3] says that the graph was discovered by Sims,  and calls it the Sims-Gewirtz graph.
\smallskip

My own interest in strongly regular graphs dates from the late 1960s, when I was told by John
McKay about the exciting discoveries of new simple groups. The topic (but not the Gewirtz graph)
is mentioned in a paper I gave at the 1969 Oxford Conference [1].
 I do not wish to claim any originality for myself, but
I am fairly sure that I initially derived my knowledge of the Gewirtz graph from a 1965
paper of W.L. Edge `On some implications of the geometry of the 21-point plane' [7].
In that paper the three sets of 56 ovals in PG$(2,4)$ are clearly described, with the
critical property that any one of the sets of 56 has the property that
two of them intersect in 0 or 2 points.

\bigskip

{\bf References}
\medskip

\begin{enumerate}

\item   N.L. Biggs. Intersection matrices for linear graphs. In: {\it Combinatorial Mathematics and
its Applications}  (ed. D.J.A. Welsh), Academic Press, 1971, 15-23.

\item   N.L. Biggs. Strongly regular graphs with no triangles. {\it arXiv}  0911.2160v1, September 2009.
Families of Parameters for SRNT Graphs. {\it arXiv} 0911.2455v1,  October 2009.

\item  A.E. Brouwer. Sims-Gewirtz graph. www.win.tue.nl/aeb/graphs/Sims-Gewirtz.html (accessed 10/11/09).

\item  A.E. Brouwer, W. Haemers. The Gewirtz Graph: an exercise in the theory of graph spectra.
{\it Europ. J. Combinatorics} 14 (1993) 397-407.

\item P.J. Cameron, J. van Lint.  {\it Designs, Graphs, Codes and their Links}. Cambridge University Press, 1991.

\item H.S.M. Coxeter. Self-dual configurations and regular graphs.  {\it Bull. Amer. Math. Soc.} 56 (1950)
413-458.

\item W.L. Edge. Some implications of the geometry of the 21-point plane. {\it Math. Zeitschr.} 87 (1965)
348-362.

\item A. Gewirtz. The uniqueness of $g(2,2,10,56)$. {\it Trans. New York Acad. Sci.} 31 (1969) 658-675.

\item A. Gewirtz. Graphs with maximal even girth. {\it Canad. J. Math} 21 (1969) 915-934.

\item C.D. Godsil. Problems in algebraic combinatorics. {\it Elect. J. Combinatorics} 2 (1995) F1.

\item C.D. Godsil, G.F. Royle. {\it Algebraic Graph Theory}. Springer, 2001.

\item  Wolfram MathWorld.  mathworld.wolfram.com/GewirtzGraph.html (accessed 10/11/09).

\end{enumerate}

\end{document}